\documentclass[10pt, twoside]{article}
\usepackage{bcc24book,url,booktabs}

\bcctitle  {Some old and new problems in combinatorial
geometry I: \\ Around Borsuk's problem}
\bccshorttitle{Around Borsuk's problem}
\bccname{Gil Kalai\footnote{Work supported by an ERC advanced grant, an ISF grant, and an NSF grant.}}
\bccshortname{G. Kalai}

\bccaddressa{Einstein Institute of Mathematics, Hebrew University of Jerusalem,\\
and Department of Mathematics, Yale University.}
\bccemaila{kalai@math.huji.ac.il}

\newtheorem{prob}[example]{Problem}

\begin{document}
\makebcctitle

\begin{abstract}

Borsuk  \cite {Bor33} asked in 1933 if every set of diameter 1 in $\R^d$ can be covered
by $d+1$ sets of smaller diameter.
In 1993, a negative solution, based on a theorem by Frankl and Wilson \cite {FW},
was given by Kahn and Kalai \cite {KK:bor}.
In this paper I will present questions related to Borsuk's problem.

\end{abstract}
\thispagestyle{empty}
\setcounter{page}{1}
\let\Horig\H


\newcommand{\bea}{\begin{eqnarray}}
\newcommand{\ba}{\begin{array}}
\newcommand{\bean}{\begin{eqnarray*}}
\newcommand{\ea}{\end{array}}
\newcommand{\eea}{\end{eqnarray}}
\newcommand{\eean}{\end{eqnarray*}}
\newcommand{\be}{\begin{equation}}
\newcommand{\ee}{\end{equation}}
\newcommand{\lb}{\linebreak}
\newcommand{\nlb}{\nolinebreak}
\newcommand{\zb}{\overline{z}}
\newcommand{\perm}{\mbox{\rm perm}}
\newcommand{\nl}{\newline}
\newcommand{\ev}{\mbox{\bf E}}
\newcommand{\lab}{\label}
\def\ct{\tilde{C}}
\def\Cox{\hfill \Box}
\def\qed{$\hfill \Box$}
\def\diseq{\, {\stackrel {{\cal D}} {=}}}
\def\dconv{\, {\stackrel {{\cal D}} {\to}}}
\def\sif{\sigma \mbox{\rm -field}}
\def\dd{\delta}
\def\ep{\epsilon}
\def\E{{\bf{E}}}
\def\Mid{\left\vert \right.}
\def\P{{\bf{P}}}
\def\N{{\bf N}}
\def\C{{\mathbb{C}}}
\def \K{{\mathcal K}}
\def\L{{\mathcal{L}}}
\def\Z{{\bf{Z}}}
\def\ZZ{{\mathbb{Z}}}
\def\F{{\mathcal{F}}}
\def\G{{\mathcal{G}}}
\def\X{{\mathcal X}}
\def\UU{{\mathbf U}}
\def\mc{\Phi}
\def\1{{\mathbf 1}}
\def\XZ{(Y,Z,\rproj)}
\def\en{\E^{\nu}}
\def\emu{\E^{\mu}}
\def\rseq{S}
\def\rproj{W}
\def\rprojf{w}
\def\rsf{{\cal G}}
\def\romega{r}
\def\rt{t}
\def\seq{{\R^\infty}}
\def\ren{\Upsilon}
\def\typren{\delta}
\def\oseq{\Omega}
\def\oproj{Y}
\def\osf{{\cal F}}
\def\oomega{\omega}
\def\measn{\alpha}
\def\measa{\beta}
\def\measb{\pi}
\def\bb{\overline{\beta}}
\def\cross{{\nearrow\mkern-18mu \searrow}}
\def\bit{\{0,1\}}
\def\toss{\{-1,1\}}
\def\x{{x_{0}}}
\def\R{{\bf R}}
\def\Ga{\Gamma}
\def\param{\varphi}
\def\whpsi{\widehat{\psi}}
\def\B{{\cal{B}}}
\def\|{\, |\! | \, }
\def\one{{\bf 1}}
\def\al{\alpha}
\def\la{\lambda}
\def\th{\theta}
\def\wh{\widehat}
\def\om{\omega}
\def\T{{\bf T}}
\def\phi{\varphi}
\def\tail{{\cal T}}
\def\and{\, \mbox{ and } \,}
\def\given{\; | \; }
\def\iid{i.i.d.\ }
\def\b2b{\overline{\beta^2}}
\def\var{{\tt Var}}
\def\cov{{\tt Cov}}
\def\res{L}
\def\shift{\Theta}
\def\M{M}
\def\coins{\{-1,1\}}
\def\bin{\{0,1\}}
\def\bias{\theta}
\def\renp{\Delta}
\def\rent{\delta}
\def\np{Y}
\def\nt{y}
\def\legt{\Lambda^{*}}
\def\sp{Z}
\def\st{z}
\def\yvars{(\np,\sp,\renp)}
\def\ma{\mu_{\measa}}
\def\HH{{\mathcal H}}
\def\mb{\mu_{\measb}}
\def\Q{{\bf Q}}
\def\omt{\tilde{\Omega}}
\def\inprob{\stackrel{\rm Pr} {\longrightarrow}}
\def\measu{\eta}
\def\rtil{\tilde{R}}
\def\muu{\mu_\eta}
\def\ls{\widehat{R}}
\def\c{\zeta}
\def\Var{{\tt Var}}

\newcommand{\cL}{{\mathcal L}}
\newcommand{\diam}{{\rm diam}}
\newcommand{\with}{\hbox{ {\rm with} }}
\renewcommand{\and}{\hbox{ {\rm and} }}
\newcommand{\comp}{\hbox{\bf CO}}
\newcommand{\Tm}{T_{{\rm mix}}}
\renewcommand{\C}{{\mathcal{C}}}
\newcommand{\p}{\mbox{\bf p}}
\newcommand{\LL}{L}
\newcommand{\EE}{{\mathcal E}}
\newcommand{\RR}{{\mathcal R}}
\newcommand{\hh}{h}
\newcommand{\lr}{\leftrightarrow}

\def\edge{{E}}
\def\STAR{\star}
\def\CYCLE{\diamondsuit}
\def\CHI{\mathchoice%
{\raise2pt\hbox{$\chi$}}%
{\raise2pt\hbox{$\chi$}}%
{\raise1.3pt\hbox{$\scriptstyle\chi$}}%
{\raise0.8pt\hbox{$\scriptscriptstyle\chi$}}}
\def\smalloplus{\raise1pt\hbox{$\,\scriptstyle \oplus\;$}}
\def\h{\widehat}
\def\H{\ell^2_-(\edge)}
\def\I#1{{\bf 1}_{#1}}

\newcommand{\Poi}{\mbox{\rm Poi}}
\newcommand{\intc}{\int_0^{2\pi}}
\def\summ{{\sum\limits}}
\def\intt{{\int\limits}}
\def\prodd{{\prod\limits}}
\def\mb{\mbox}
\def\l{\left}
\def\r{\right}
\def\sig{\sigma}
\def\lam{\lambda}
\def\alp{\alpha}
\def\eps{\epsilon}
\def\tends{\rightarrow}
\def\X{{\mathcal X}}

\newcommand{\cG}{{\mathcal{G}}}
\newcommand{\cB}{{\mathcal{B}}}
\newcommand{\GPC}{\cG_{\textsc{pc}}}
\newcommand{\cF}{{\mathcal{F}}}
\newcommand{\A}{{\mathcal A}}
\newcommand{\oneb}[1]{\one_{\{ #1 \}}}
\newcommand{\deq}{\stackrel{\scriptscriptstyle\triangle}{=}}
\newcommand{\GC}{{\mathcal{C}_1}} 
\newcommand{\tGC}{{\tilde{\mathcal{C}}_1}}
\newcommand{\TC}[1][\mathcal{C}_1]{#1^{(2)}} 
\newcommand{\hGC}{{\widehat{\mathcal{C}}_1}}
\newcommand{\Po}{\operatorname{Po}}
\newcommand{\Bin}{\operatorname{Bin}}
\newcommand{\Geom}{\operatorname{Geom}}











\section {Introduction}

The title of this paper is borrowed from Paul Erd\Horig{o}s who used it
(or a similar title) in many lectures and papers, e.g., \cite {E2}.
I will describe several open problems in the interface between combinatorics and geometry, mainly convex geometry.
In this part, I describe and pose questions related to the Borsuk conjecture. The selection of problems
is based on my own
idiosyncratic tastes. 
For a fuller picture, the reader is
advised to read the review papers on Borsuk's problem and related questions by Raigorodskii
\cite {Rai07,Rai14,Rai13,Rai04,Rai07b}. Among other excellent sources are \cite {BG,BMS,BMP,Mat:b,Pak:b}.

Karol Borsuk \cite {Bor33} asked in 1933 if every set of diameter 1 in $\R^d$ can be covered
by $d+1$ sets of smaller diameter. That the answer is positive
 was widely believed and referred to as the \emph{Borsuk conjecture}.
However, some people, including Rogers, Danzer, and Erd\Horig{o}s,
suggested that a counterexample might be obtained from some clever combinatorial configuration.
In hindsight, the problem is related to several
questions that Erd\Horig{o}s asked
and its solution was a great triumph for Erd\Horig{o}sian
mathematics. 

\section{Better lower bounds to Borsuk's problem}
\subsection {The asymptotics}

Let $f(d)$ be the smallest integer such that every set of diameter
one in $\R^d$ can be covered by $f(d)$ sets of smaller diameter.
The set of vertices of a regular simplex of
diameter one demonstrates that
$f(d) \ge d+1$. The famous Borsuk--Ulam theorem \cite {Bor33} asserts that the $d$-dimensional ball of diameter 1 cannot be
covered by $d$ sets of smaller diameter. The Borsuk--Ulam theorem has many important applications in many areas of mathematics.
See Matousek's book \cite {Mat} for applications and connections to combinatorics.
In the same paper \cite {Bor33} Borsuk asked if $f(d) \le d+1$. This was proved for $d=2,3$.
It was shown by Kahn and Kalai
 \cite {KK:bor} that  $f(d) \ge 1.2^{\sqrt d}$,
by Lassak
\cite {Las} that $f(d) \le 2^{d-1}+1$ and by Schramm
\cite {Sch88a} that $f(d) \le (\sqrt{3/2}+o(1))^d$.

\begin {prob}
Is $f(d)$ exponential in $d$?
\end {prob}

The best shot (in my opinion) at an example leading to a positive answer is:

\begin {itemize}

\item[(a)] Start with binary linear codes of length $n$ (based on algebraic-geometry codes)
with the property that the number of maximal-weight codewords is exponential in $n$.

\item[(b)] Show that the code cannot be covered by
less than an exponential number of sets that do not realize the maximum distance.

\end {itemize}

Part (a) should not be difficult, given that it is known that
for certain AG-codes the number of minimal-weight codewords is exponential in $n$ \cite {ABV}.
Part (b) can be difficult, but the algebraic techniques used for the Frankl and Wilson
theorem may apply.


The Kahn--Kalai counterexample and many of the subsequent results
depend on the Frankl--Wilson \cite {FW} theorem or on some
related algebraically-based combinatorial results.
(One can rely also on the Frankl--R\"odl theorem \cite {FR},
which allows much greater generality but not as good quantitative estimates.)
We will come back to these results later on.

Let $g(d)$ be the smallest integer such that every {\it finite} set of diameter
one in $\R^d$ can be covered by $g(d)$ sets of smaller diameter.

\begin {prob}
Is $f(d) = g(d)$?
\end {prob}

I am not aware of any reduction from infinite sets to finite sets,
and indeed the proof of Borsuk's conjecture for $d=2,3$ is easier if one
considers only finite sets. On the other hand,
the counterexamples are based on finite configurations.
Perhaps one can demonstrate a gap between the finite and infinite behavior for
some extension or variation of the problem, e.g., for arbitrary metric spaces.
(Our knowledge of $f(d)$ does not seem accurate enough to hope to
prove that such a gap exists for the original problem.)

\subsection {Larman's conjecture}

The counterexample to Borsuk's conjecture is based on the  special case where the set consists of 0-1 vectors
of fixed weight. Here, the conjecture has an appealing combinatorial formulation.

\begin {prob}[Larman's conjecture]
Let $\cal F$ be a $t$-intersecting family of $k$ sets from $[n]$. Then $\cal F$ can be covered
by $n$ $(t+1)$-intersecting subfamilies.
\end {prob}

We now know that Larman's conjecture does not hold in general. However:

\begin {prob}
Is Larman's conjecture true for $t=1$?
\end {prob}

For more discussion of the combinatorics of Larman's conjecture and
related combinatorial questions on the packing and coloring of graphs and hypergraphs,
see \cite{Kah94}.
We can sort of ``dualize''
the $t=1$ case of Larman's conjecture by replacing ``intersecting''
(i.e., ``every pair of sets has at least one common element'')
by ``nearly disjoint'' (namely,  ``every pair of distinct
sets has at most one common elements'') and thus recover the famous:

\begin {conj} [Erd\Horig{o}s--Faber--Lov\'asz]
Let $\cal F$ be a  family of nearly disjoint $k$-sets from $[n]$.
Then $\cal F$ is the union of $n$ matchings.
\end {conj}

While the Erd\Horig{o}s--Faber--Lov\'asz conjecture is still open it is known that $n$
is the right number for the
fractional version of the problem \cite {KS},
and
that $(1+o(1))n$ matchings suffice \cite{Kah92}.
Such results are not available for the $t=1$ case of Larman's conjecture but a counterexample to
a certain strong form of the conjecture is known \cite {KK:s}.


\subsection {Embedding the elliptic metric into an Euclidean one}

Let ${\cal E}_n$ be the {\em elliptic} metric space of lines through the origin in $\R^n$ where the distance between
two lines is the (smallest) angle between them.
So the diameter of ${\cal E}_n$ is $\pi/2$ and
the famous Frankl--Wilson theorem implies that ${\cal E}_n$ cannot be covered by less
than an exponential number (in $d$) of
sets of diameter smaller than $\pi/2$. The proof by Kahn and Kalai can be seen as adding a single simple fact:
$\cal E$ can be  embedded into an Euclidean space $\R^{n(n+1)/2}$ by the map $$x \to x \otimes x.$$
The distance between $x \otimes x$ and $y \otimes y$ is a simple
monotonic function of the distance between $x$ and $y$.
(Here $\|x\|_2=1$ and we note that $x$ and $-x$ are mapped onto the same point.)

The original counterexample, $C_1$, is the image under this map
for (normalized) $\pm 1 $ vectors of length $n$ with $n/2$ '1's ($n$ divisible by 4).
Another example, $C_2$, is the image of all $\pm 1$ vectors, and
we can also look at $C_3$ the image of all unit vectors.
All these geometric objects are familiar: $C_3$ is the unit vectors
in the cone of rank-one positive semi-definite matrices, $C_2$ is called the {\em cut polytope,}
and $C_1$ is the polytope of balanced cuts \cite {DL}.

We can ask if there are more economic embeddings
of the elliptic space into a Euclidean space. Namely,
is there an embedding $\phi: {\cal E}_n \to \R^m$, $m=o(n^2)$, such that
$\|\phi (x) - \phi(y)\|_2= \phi (d(x,y))$ for some (strictly) monotone function $\phi$?

The answer to this question is negative by an important theorem of de Caen
from 2000 \cite {dC}.

\begin {theorem}[de Caen]
There are quadratically many equiangular lines in ${\cal E}_n$.
\end {theorem}

Weaker forms of embeddings of ${\cal E}_n$ into Euclidean spaces possibly with some symmetry-breaking
may still lead
to improved lower bounds for $f(d)$, and are of independent interest.


\begin {prob}
Is there a continuous map $\phi: {\cal E}_n \to \R^m$, $m=o(n^2)$ so that $\phi$
preserves the set of diameters of ${\cal E}_n$?
\end {prob}

\subsection{Spherical sets without pairwise orthogonal vectors}

Regarding ${\cal E}_n$ itself,
Witsenhausen
asked in 1974 \cite {Wit74} what is the maximum volume $\mu (A)$ of an $n$-dimensional
spherical set $A$
without a pair of orthogonal vectors. Witsenhausen proved that:

$$\mu (A) \le 1/(n+1).$$

The following natural conjecture is very interesting:

\begin {conj}
\label {c:caps}
Let $A$ be a measurable subset of $S^n$ and suppose that $A$ does not contain two orthogonal vectors.
Then the volume of $A$ is at most twice the volume of two spherical caps of radius $\pi/4$.
\end {conj}

Asymptotically this conjecture asserts that a subset of the $n$ sphere of measure $(1/\sqrt 2+o(1))^n$ must contain a pair
of orthogonal vectors. If true, this can replace the Frankl--Wilson bound and
will show that $C_3$ defined above is  a
counterexamples to Borsuk's conjecture for $d>70$ or so.
The Frankl--Wilson theorem gives that if $\mu(A) > 1.203...^{-n}$ then $A$ contains two orthogonal vectors.
It seems that the main challenge is to extend the
linear algebra/polynomial method from 0-1 vectors to general vectors.
One important step was taken by Raigorodskii \cite {Rai99} who improved
the bound to $1.225^{-n}$.

Remarkably, the upper bound of 1/3 for the two-dimensional case stood unimproved for 40 years until very
recently DeCorte
and Pikhurko improved it to 0.31.. --\cite {DP14}! The proof uses
Delsartes' linear programming method \cite {delsarte73} combined with a combinatorial argument.

\subsection {Borsuk's problem for spherical sets}

Borsuk's problem itself has an important extension to spherical sets.
Consider a set of Euclidean diameter 1 on a $d$-dimensional sphere $S^{d-1}_r$ of radius $r$.

\begin {prob}
What is the maximum number $f_r(d)$ of parts one needs to
partition any set of diameter 1 on $S_r^d$?
\end {prob}

Obviously, one has $f_r(d)\le f(d)$ for any $r$, and we as well have $f_{1/2}(d)=d+1$ due to Borsuk--Ulam theorem.

Kupavskii and Raigorodskii \cite {KR:ball}  proved the following theorem:

\begin{theorem}

Given $k\in \N$, if $r>\frac12\sqrt{\frac{2k+1}{2k}}$, then there exists $c>1$
such that $f_r(d)\ge(c + o(1))^{\sqrt[2k]d}$. Moreover, there exist a $c>0$ such that if $r>1/2+c\frac{\log\log d}{\log d}$,
then for all sufficiently large $d$ we have $f_r(d)>d+1$.

\end{theorem}

The proof is based on mappings involving multiple tensor products.
We note that
embeddings of similar nature
via multiple tensor-products play a role also in the disproof of Khot and Vishnoi \cite {KV}
of the Goemans--Linial conjecture.

\begin {prob}
Is it the case that for every $r<1/2$ there is a constant $C_r>1$ such that
$$f_r(d)\ge C_r^d?$$
\end {prob}

\subsection {Low dimensions and two-distance sets}

The initial counterexample showed that the Borsuk conjecture
is false for $n=1325$ and all $n>2014$ and there were gradual
improvements over the years down to 946 (Nilli \cite {N}),
903 (Weissbach), 561 (Raigorodski \cite {Rai245}), 560 (Weissbach \cite {Wei}),
323 (Hinrichs \cite {Hin}), 321 (Pikhurko \cite {Pikh}), 298 (Hinrichs and Richter \cite {HC}).
The construction of Hinrichs and the subsequent ones remarkably rely on the Leech lattice.

A two-distance set is a set of vectors in $\R^d$ that attain only two distances.
Larman asked early on 
and asked again recently:

\begin {prob}
Is the Borsuk conjecture correct for two-distance sets?
\end {prob}

This has proven to be a very fruitful question. In 2013 Bondarenko \cite {Bon}
found a two-distance set
with 416 points in 65 dimensions that cannot be
partitioned into less than 83 parts of smaller diameter. Remarkable!
Jenrich \cite {Jen} pushed the dimension down to 64. These constructions beatifully relies
on known strongly-regular graphs.

\begin {prob}
What is the smallest dimension for which Borsuk's conjecture fails? Is Borsuk's conjecture correct
in dimension 4?
\end {prob}

In dimensions 2 and 3 Borsuk's conjecture is correct. Eggleston gave the first proof for dimension 3 \cite {Egg},
which was followed by simpler proofs  by Gr\"unbaum \cite {Gr2} and Heppes \cite {Hep}.
A simple proof for finite sets of points in 3-space was found by Heppes and Revesz \cite {HR}.
For dimension 2 it follows from an earlier 1906 result that every set of diameter one can be embedded into a regular
hexagon whose opposite edges are distance one apart.
For a simpler argument see  Pak's book \cite {Pak:b}. Here too, for finite
configurations the proof is very simple.

\section {Upper bounds for Borsuk's problem and sets of constant width}

\subsection {Improving the upper bound}

Lassak \cite {Las} proved that for every $d$, $f(d) \le 2^{d-1}+1$ (and this still gives
the best-known bound when
the dimension is not too large).
Schramm \cite {Sch88a} proved that every convex body of constant width 1 can be covered by
$(\sqrt{3/2}+o(1))^d$ smaller homothets. It is a well-known fact \cite {Egg:b}
that every set of diameter
one is contained in a set of constant width 1, and, therefore, for proving an
upper bound on $f(d)$ it is enough to
consider sets of constant width.
Bourgain and Lindenstrauss \cite {BL} showed that
every convex body in $\R^d$ of diameter 1 can be covered by $(\sqrt{3/2}+o(1))^d$ balls of diameter 1.
Both these results show that $f(d) \le (\sqrt{3/2}+o(1))^d$.

\begin {prob}
Prove that
$f(d) \le C^d$ for some $C < \sqrt{3/2}$.
\end {prob}

We note that Danzer constructed a set of diameter 1  that  requires exponentially many balls to cover.
Danzer constructed a set for which  $(1.003)^d$ balls needed, and Bourgain and Lindenstrauss in 1991 \cite {BL}
found much better bound,  $(1.064)^n$.

As for covering by smaller homothets we recall the famous:

\begin {conj}[Hadwiger \cite {Had1,Had2}] Every convex body $K$ in $ \R^d$ can be covered by $2^d$ smaller homothets of $K$.
\end {conj}

The case of sets of constant width is of particular interest:

\begin {prob} Are there $\epsilon>0$ and sets of constant width in $\R^n$ that
require at least $(1+\epsilon)^n$ smaller homothets to cover?
\end {prob}

Note that a positive answer neither implies nor follows from a $(1+\epsilon)^d$ lower bound for the Borsuk number $f(d)$.

\subsection{Volumes of sets of constant width}

Let us denote the volume of the $n$-ball of radius 1/2 by $V_n$.

\begin {prob}[Schramm \cite{Sch88b}] Is there
some $\epsilon >0$ such that for every $d>1$
there exists a set $K_d$ of constant width 1 in dimension $d$
whose volume satisfies $VOL(K_d) \le (1-\epsilon)^d V_d$?
\end {prob}

Schramm raised a similar question for spherical sets of constant width and pointed out that a
negative answer for spherical sets will push
the $(3/2)^{d/2}$ upper bound for $f(d)$ to $(4/3)^{d/2}$.

\section {Saving the Borsuk conjecture}

\subsection {Borsuk's conjecture under transversality}

I would like to examine the possibility that Borsuk's conjecture is correct except for
some ``coincidental'' sets.
The question is how to properly define ``coincidental,'' and we will now give it a try!

Let $K$ be a set of points in $\R^d$ and let $A$ be a set of pairs of points in $K$. We say that the pair $(K, A)$
is {\em general} if for every continuous deformation of the distances on $A$ there is a deformation $K����$ of $K$
which realizes the deformed distances.

\begin{rem} This condition is related to the ``strong Arnold property'' (a.k.a. ``transversality'')
in
Colin de Verdi\'ere's theory of  invariants of graphs \cite {CdV}.\end{rem}

\begin {conj}
If $D$ is the set of diameters in $K$ and $(K,D)$ is general then $K$ can be partitioned into $d+1$ sets of smaller diameter.
\end {conj}

We further propose (somewhat more strongly) that this conjecture
holds even when ``continuous deformation'' is replaced
with ``infinitesimal deformation.''

The finite case is of special interest. A graph embedded in $\R^d$ is {\it stress-free}
if we cannot assign not-all-zero weights to the edges such that the weighted
sum of the edges containing any  vertex $v$ (regarded as vectors from $v$) is zero for every vertex $v$. Here we embed the vertices and
regard the edges as straight line segments. (Edges may intersect.) Such a graph is called a ``geometric graph.''
When we restrict the conjecture to finite configurations of points we get:

\begin {conj}
If $G$ is a stress-free geometric graph of diameters in $\R^d$  then $G$ is $(d+1)$-colorable.
\end {conj}

\begin {rem} A stress-free graph for embeddings into $\R^d$ has at most
$dn-{{d+1} \choose {2}}$ edges and therefore its chromatic number is at most $2d-1$.
\end {rem}



\subsection {A weak form of Borsuk's conjecture}

\begin {conj}
Every polytope $P$ with $m$ facets can be covered by $m$ sets of smaller diameter.
\end {conj}

This conjecture was motivated by recent important works on projections of polytopes \cite {FMPTdW}. A positive answer will give an alternative path
for showing that the cut polytope cannot be described as a projection of a polytope with only polynomially many facets.

\subsection {Classes of bodies for which Borsuk's conjecture holds}

Perhaps the most natural way to ``save'' Borsuk's conjecture is given by:

\begin {prob}
Find large and interesting classes of convex bodies for which Borsuk's conjecture holds!
\end {prob}

Borsuk's conjecture is known to be true for centrally symmetric bodies,
Hadwiger proved it for smooth convex bodies \cite {Had2}, and
Boris Dekster proved the conjecture both for bodies of revolution \cite {dec2} and for
convex bodies with a belt of regular points \cite {dec3}.


\subsection {Partitioning the unit ball and diametric codes}

The unit ball in $\R^d$ can be covered by $d+1$ convex sets of smaller diameter. But how much smaller?
We do not know the answer. Let $u(d)$ be the minimum value of $t$ such that the unit ball
in $\R^d$ can be covered by $d+1$ sets of diameter at most $t$.

\begin {prob}
Determine the behavior of $u(d)$!
\end {prob}

The motivation for this question comes from an even stronger form of Borsuk's conjecture asserting that
every set of diameter 1 can be covered by $d+1$ sets of diameter $u(d)$. It was also conjectured that
the optimal covering for the sphere is described by
a partition based on the Voronoi regions of a regular simplex that gives $u(n) \le 1-\Omega (1) /n. $
This is known to be optimal in dimensions two and three and is open in higher dimensions.
Larman and Tamvakis \cite {LT} showed by a volume argument
that $u(n) \ge 1-3/2 \log n/n +O(1/n)$. See also \cite {dec1}.

It will be interesting to close the logarithmic gap for $u(d)$.
I don't know what one should expect for the answer, and it will be quite exciting
if the standard example is {\em not} optimal.

We can more pose general questions:

\begin {prob}

(i) What is the smallest number of sets of diameter $t$ that are needed to cover the unit sphere?

(ii) What is the largest number of convex sets of width $\ge t$ that can be packed into the unit sphere? (The
\emph{width} of a convex set is the minimum
distance between opposite supporting hyperplanes.)

\end {prob}

Let $\Omega_n=\{0,1\}^n$. We can ask the analogous questions about the binary cube.

\begin {prob}
\label{4.7}

(i) What is the smallest number of sets of diameter $t$ that are needed to cover $\Omega_n$?

(ii) What is the largest number of sets of width $\ge t$ that can be packed into $\Omega_n$?

\end {prob}

Here by ``diameter'' and ``width'' we refer to the Euclidean
notions (for which, for Problem \ref {4.7}, ``diameter'' essentially coincides with the
Hamming diameter).

\section {Unit-distance graphs and complexes}

For a subset $A$ of $\R^n$ the {\em unit-distance graph} is a
graph whose set of edges consists of pairs of points of $A$ of distance 1.
If all pairwise distances are at most 1, we call the unit-distance graph a {\em diameter graph}.
If all the pairwise distances are at least 1 we call it a {\em  kissing graph}.
Borsuk's question is a question about coloring diameter graphs.

\begin {prob}
\label {p:distancegraph}
What is the maximum number of edges, the maximum chromatic number, and the maximum minimal degree
for the diameter graph, kissing graph,
and unit-distance graph for a set of $n$ points in $\R^d$?
\end {prob}

Finding the maximum number of edges in a planar unit distance graphs
is a famous problem by Erd\Horig{o}s \cite {E}. Another famous problem by Hadwiger and Nelson
is about the chromatic number of the planar unit distance graph and yet another famous question is if
the minimal kissing number of a set of $n$ points in $\R^d$ is exponential in $d$, see \cite {conway93,A:kissing}.

We can define also the {\em unit-distance complex}
to be the simplicial complex of cliques in the unit-distance graphs or, alternatively,
the simplicial complex whose faces are sets of points in $A$ that form regular simplices of diameter 1. And again
when the diameter of $A$ is 1 we call it the {\it diameter complex} and when the minimum distance is 1 we call it the
{\it kissing complex.}



\begin {prob}
What is the maximum number of $r$-faces for the diameter complex, kissing complex,
and the unit-distance complex for a set of $n$ points in $\R^d$?
\end {prob}

For the chromatic number of the unit-distance graph it makes
a difference if we demand further that each color class be
measurable. (This is referred to as the {\it measurable chromatic number}.)
For progress on the chromatic number
of unit-distance graphs, see \cite {Rai11,Kup1,Kup2,Kup3,KR}. For progress on the
measurable chromatic number and related questions, see
\cite {bachoc09,bachoc13,bachoc14,delaat13}.

Rosenfeld asked (see \cite {Ros}):

\begin {prob}
Does the graph whose vertex set is the set of points in the plane and whose edges represent points whose distance is an odd integer
have a bounded chromatic number?
\end {prob}

For measurable chromatic numbers the answer is negative as follows from a
theorem of Furstenberg, Katznelson and Weiss 
that asserts that every planar set of positive measure realizes all sufficiently large distances. See also
\cite {Ste} for a simple direct proof.

\subsection {Schur's conjecture}

A conjecture by Schur deals with an interesting special case:

\begin {conj} [Schur]
The number of $(d-1)$-faces of every diameter complex for a set of $n$ points in $\R^d$ is at most $n$.
\end {conj}

The planar case is an old result
and it implies
a positive answer to Borsuk's problem for finite planar sets.
The proof is based on an observation that sets the metric aside:
the edges of the diameter graph are pairwise intersecting and therefore we need to show that
every geometric graph with $n$ vertices and $n+1$ edges
must have two disjoint edges. This result by Hopf and Pannwitz \cite {HP} from 1934 can be
seen as the starting
point of ``geometric graph theory'' \cite {Pach}.
Zvi Schur was a high school teacher who did research in his spare time.
He managed to prove his conjecture in dimension 3
(see \cite {SKMP}) but in his writing he mentioned that ``the power of my methods diminishes as the dimension goes up.''
The paper \cite {SKMP} includes also a proof that in any dimension
the number of $d$-faces of the diameter complex is at most one.
Schur's  conjecture has recently been proven
by Kupavskii and Polyanski \cite{KP}!
(For $d=4$ it was proved by Bulankina, Kupavskii, and Polyanskii \cite {BKP}.)
A key step in Kupavskii and Polyanskii's work is proving the  $k=m=d-1$ case of the following additional
conjecture by Schur, still open in the general case.

\begin {conj} [Schur]
Let $S_1$ and $S_2$ be two regular simplices of dimensions $k$ and $m$ in $\R^d$
such that their union has diameter 1. Then
$S_1$ and $S_2$  share at least $\min (0,k+2m-2d+1)$ vertices for $k \ge m$.
\end {conj}

Heppes and R\'ev\'esz proved that the number of edges in the diameter graph of $n$ points in space is $2n-2$.
This gives an easy proof of Borsuk's conjecture for finite sets of points in $\R^3$.


A natural weakening of Borsuk's conjecture is:

\begin {prob}
What is the smallest $r=r(d)$ such that every set of diameter $1$ in $\R^d$ can be covered by $d+1$ sets,
none of which contains an $r$-dimensional simplex of diameter 1?
\end {prob}

Unit-distance graphs and especially diameter
graphs and complexes are closely related to the study of {\it ball polytopes}.
Those are
convex bodies that can be described as the intersection of unit balls.
A systematic study of ball polytopes was initiated
by C\'aroly Bezdek around 2004 and they were also studied by Kupitz, Martini,
and Perles, see \cite {BLN,BN,KMP}. Ball polytopes
are also related to sets of constant width.

\subsection {Tangent graphs and complexes for collections of balls (of different radii)}

We can try further to adjust the problems discussed in this section to the case where we have a collection
$A$ of points in $\R^d$
and a close ball centered around each point. Two balls can be in three mutual positions (that we care about):
They can be disjoint, they can have intersecting interiors, or they can be tangential.

The {\em tangent graph} is a graph whose set of vertices is $A$ and a
pair of vertices are adjacent if the corresponding balls are tangential.
Note that if all balls have the same radius 1/2, then the tangent graph is the unit-distance graph.
As before we can consider also the {\it tangent complex} - the simplicial complex described by
cliques in the tangent graph.

\begin {prob}

(i) What is the maximum number of edges in a tangent graph (especially in the plane)?
What is its maximum chromatic number (especially in the plane)?

(ii) If every two balls intersect, then the tangent graph is a generalization of the diameter graph.
Again we can ask for the maximum number of edges, cliques of size $r$, and the chromatic number. Again we can
ask if when the graph is stress-free the chromatic number is at most $d+1$.

(iii) If every two balls have disjoint interiors then the tangent graph is a generalization of the kissing graph.
Again we can ask for the maximum number of edges, cliques of size $r$, the maximum minimal degree, and
the chromatic number.
\end {prob}

In the plane we can find $n$ points and $n$ lines with $n^{C4/3}$ incidences and the famous
Szemeredi--Trotter theorem (see, e.g., \cite {Szek})
asserts that this is best possible. Now, we can replace each point by a small circle,
arrange for the lines incident to the points to be tangential to them, and regard the lines as circles as well.
This shows
that  tangent graphs with $n$ vertices in the plane can have as much as $n^{C4/3}$ edges.
It is conjectured by Pinchasi, Sharir,
and others that

\begin {conj}
(i) Planar tangent graphs with $n$ vertices can have at most\\ $n^{4/3} polylog (n)$ edges.

(ii)  More generally, $m$ red discs and $n$ blue discs (special case: $n$ blue points),  can touch
at most $((mn)^{2/3}+m+n) polylog (m,n)$ times.
\end {conj}

This conjecture proposes a profound extension of the Szemeredi--Trotter theorem.
The best known upper bound $n^{3/2}\log n$ is by Markus and Tardos \cite {MT}, following an
earlier argument by Pinch\'asi and Radoici\'c \cite {PR03}.
This particular approach based on a certain ``forbidden
configurations'' -- a self crossing 4-cycle -- cannot lead to better exponents.
Sharir found a beautiful connection with Erdos's distinct distances problem \cite {GK}
which also shows that the
$polylog (n,m)$ term cannot be eliminated. Indeed, assume you have $n$ points with just $x$ distances.
Then draw around each point $x$ circles whose radii are the $x$ possible distances and then you get a collection of $m=nx$ circles 
and $n$ points with $n^2$ incidences (because every point lies on $n$ circles exactly).
Therefore: $n^2 \leq polylog (xn) (n(xn))^2/3+n+xn$ which implies $x \geq n polylog (n)$.

For circles that pairwise intersect, Pinchasi \cite {P02} proved a Gallai--Sylves\-ter
conjecture by Bezdek asserting that (for more than 5 circles)
there is always a circle tangential to at most
two other circles. This
was the starting point of important studies \cite {ALPS,ANPPSS} concerning arrangement of
circles and pseudo-circles in the plane. Alon, Last, Pinchasi, and Sharir \cite {ALPS} showed
 an upper bound of $2n-2$ for the number of edges in the tangent graph for pairwise intersecting circles.

The problem considered in this section can be asked under greater generality in at least two ways:
one important generalization is to consider two circles adjacent if their intersection is an {\em empty lens,}
that is, not intersected by the boundary of another disc.
Another generalization is for pseudocircles (where both notions of adjacency essentially coincide).

Let me end with the following problem:

\begin {conj}[Ringle circle problem]
Tangent graphs for finite collections of circles in the plane such that no more than two circles pass through a point
have bounded chromatic numbers.
\end {conj}

Recently, Pinchasi proved that without the assumption that no more than two circles pass through a point, the chromatic number
is $O(\log^2n)$, where $n$ is the number of vertices of the graphs. Pinchasi also gave an example 
where (again, dropping the extra assumption) you need $\log n$ colors.

\section {Other metric spaces}

\subsection {Very symmetric spaces}

We already discussed Borsuk's problem for spherical sets. We can also ask

\begin {prob}

Study the Borsuk problem, and other questions considered above, for very symmetric spaces
like the hyperbolic space, the Grassmanian, and $GL(n)$.

\end {prob}

The Grassmanian, the space of $k$-dimensional linear spaces of $ \R^n$ is of special interest. The ``distance'' between two vector spaces can be seen as a
vector of $k$ angles, and there may be several interesting ways to extend the questions considered here.  (The case $k=1$ brings us back to the Elliptic space).

\subsection { Normed spaces}

Given a metric space $X$ and a real number $t$ we can consider the \emph{Borsuk number} $b(X,t)$
defined as the smallest integer such that every subset of
diameter $t$ in $X$ can be covered by $b(X,t)$ sets of smaller diameter. There are interesting results and questions
regarding Borsuk's numbers of various metric spaces. Let $a(X,t)$ be the maximum cardinality of an {\em equilateral}
subset $Y \subset X$ of diameter $t$ (namely, a set so that every pairwise distance between distinct points in $Y$ is $t$).
Of course, $a(X,t) \le b(X,t)$. Understanding $a(X,t)$ for various metric spaces is of great interest. Kusner conjectured that an equilateral
set in $\ell_1^n$  has at most $2n$ elements and
an equilateral set in $\ell_p^n$  has size at most $n+1$ for $p, 1<p<\infty$. Smyth found the first
polynomial upper bound for the size of an equilateral set in $\ell_1$ which followed by
an important result by Alon and Pudlak \cite {AP03}:

\begin {theorem}[Alon and Pudlak]
For an odd integer $p$, an equilateral set in $\ell_p^n$
has at most $c_pn \log n$ points.
\end {theorem}

When we move to general normed spaces there are very basic things we do not know. It is widely conjectured that:

\begin {conj}
\label {c:ens}
Every normed $n$-dimensional space has an equilateral set of $n+1$ points.
\end {conj}

For more on this conjecture see    Swanepoel   \cite {Swa}.

Petty \cite {Pet} proved the $n=3$ case and his proof is
based on the topological fact that a Jordan curve in the plane enclosing the origin
cannot be contracted without passing through the origin at some stage.
Make\'ev proved the four-dimensional case
using more topology. Brass and Dekster proved independently a $(\log n)^{1/3}$ lower bound and a
major improvement by Swanepoel and Villa \cite {SV} improved the lower bound to
$\exp (c\log n)^{1/2}$. I would not be surprised if Conjecture \ref {c:ens} is false.
It is known that $2^n$ is an upper bound for the size of an equilateral set for a normed $n$-dimensional space.

Let me end this section with a beautiful result of Matou${\rm \check{s}}$ek
about unit distances in normed space.
One of the most famous problems in geometry is Erd\Horig{o}s' unit distance problem of finding
the maximum number of unit distances among $n$ points in the plane. This question can be asked with respect to every
planar norm with unit ball $K$. It is known that for every norm the number of edges can be as
large as $\theta(n \log n)$ and here we state a breakthrough theorem by
Matou${\rm \check{s}}$ek \cite {Mat:UD}:

\begin {theorem}
There are norms (in fact, for most norms in a Baire category sense) for which the maximum number of unit distances
on $n$ points is $O(n \log n\log \log n$).
\end {theorem}



\section {Around Frankl--Wilson and Frankl--R\"odl}

\subsection {The combinatorics of cocycles and Tur\'an numbers}

The original counterexamples to  Larman's conjecture (and Borsuk's conjecture) were based on cuts:
we consider the family of edges of complete bipartite graphs with $4n$ vertices. (In one variant we consider
balanced bipartite graphs, and in another, arbitrary bipartite graphs.) We now consider high-dimensional
generalization of cuts in graphs.

A ($(k-1)$-dimensional) {\it cocycle} is a $k$-uniform hypergraph $G$ such that every $k+1$ vertices contains an even number of edges.
Equivalently, you can start with an arbitrary $(k-1)$ uniform hypergraph $H$ and consider the $k$-uniform hypergraph $G$ of all $k$-sets
that contain an odd number of edges from $H$. Cocyles are familiar objects from simplicial cohomology and they have
also been studied by combinatorialists and mainly by
Seidel \cite {Sei}.

For even $k$, let $f(n,k)$ be the largest number of edges in a $(k-1)$-dimensional cocycle with $n$ vertices.
(Note that when $k$ is odd, the complete $k$-uniform hypergraph is a cocycle.)
Let $T(n,k,k+1)$ be the maximum number of edges in a $k$-uniform hypergraph without
having a complete sub-hypergaph with $(k+1)$ vertices.

\begin{conj}[\cite{Ka:t}] When $k$ is even, $T(n,k,k+1)=f(n,k)$.
\end {conj}

The best constructions for Tur\'an numbers $T(n,2k,2k+1)$ are obtained by cocycles.
Let me just consider the case  where $k=2$.
For a while the best example was based on a planar drawing of $K_n$
with the minimum number of crossings. For every such drawing
the set of 4-sets of points without a crossing
is an example for Tur\'an's (5,4) problem because $K_5$ is non-planar.
It is easy to see
that this non-crossing hypergraph is also a cocycle. In 1988
de Caen, Kreher, and Wiseman \cite {dCKW} found a better, beautiful example:
consider a $n/2$ by $n/2$ matrix $M$
with $\pm 1$ entries.
Your hypergraph vertices will correspond to rows and columns of $M$.
It will include all 4-tuples with 3 rows or with 3 columns
and also all sets with 2 rows and 2 columns such that the product
of the four matrix entries is -1. The expected number of edges in the hypergraph for a random $\pm 1$
matrix is $(11/16 +o(1)) {{n} \choose {4}}$.

As for upper bounds, the best-known upper bounds are stronger for cocycles.  Peled \cite {Peled}
used a flag-algebras technique to show that
$f(n,4) \le {{n} \choose {4}} (0.6916+o(1))$.

\subsection {High-dimensional versions of the cut cone and the cone of rank-one PSD matrices}

The counterexamples for Borsuk's conjecture were very familiar geometric objects
\cite {DL}.
The example based on bipartite graphs
(where the number of edges is arbitrary) is the {\em cut-polytope}.
The image of the elliptic space under the map $x \to x \otimes x$ is simply the set of unit
vectors in the cone of rank-one positive semidefinite matrices.
The unit vectors in the cone of cocycles
is an interesting generalization
of the cut polytope since for graphs (1-dimensional complexes) it gives us the cut-polytope.


\begin {prob}
Find and study a ``high-dimensional'' extension of the cone of rank one
PSD matrices (analogous to the cone of cocycles).
\end {prob}

One possibility is the following: start with an arbitrary real-valued function $g$ on ${{[n]} \choose {k-1}}$
and derive a real-valued function on ${{[n]} \choose {k}}$ by:
$$f(T)= \prod \{ g(S): S \subset T, |S|=k-1 \}.$$
Let $U_{k,n}$ be the cone of all such $g$'s.

Speculative application to Borsuk's problem is given by:

\begin {conj}
\label {c:cycles}

(i)
The set of unit vectors in the cone of 3-cocycles with $n$ vertices
demonstrates a Euclidean set in $\R^d$ that cannot be covered by less
than $exp (d^{4/5})$ sets of smaller diameter.

(ii) The set of norm-1 vectors in $U_{4,n}$ demonstrates a
Euclidean set in $\R^d$ that cannot be covered by less
than $exp (d^{4/5})$ sets of smaller diameter.

\end {conj}

\subsection {The Frankl--Wilson and Frankl--R\"odl theorems}

We conclude this paper with the major technical tool needed for the disproof of Borsuk's conjecture, which is the
Frankl--Wilson (or Frankl--R\"odl) forbidden intersection theorem. Most of the
counterexamples to Borsuk's conjecture in
low dimensions are based on algebraic techniques
``the polynomial method'' (or some variant) which seem related to the
technique used for the proof of Frankl--Wilson's theorem. (The only exception are the new
examples based on strongly regular graphs.)
The Frankl--Wilson theorem \cite {FW} is wonderful and miraculous and the Frankl--R\"odl
theorem \cite {FR}  is great
- it allows many extensions (but not with sharp constants).
The proof of Frankl--Wilson is a terrific demonstration of the linear-algebra method. The proof of
Frankl--R\"odl is
an ingenious application (bootstrapping of a kind) of isoperimetric results.
Recently Keevash and Long \cite {KL:fr} found a new proof of Frankl--R\"odl's theorem based on the Frankl--Wilson
theorem.

\begin {prob}
Is there a proof of Frankl--R\"odl's theorem based on Delsarte's linear-programming method \cite {delsarte73}?
\end {prob}

The work of Evan and Pikhurko \cite {DP14} mentioned above suggests that applying the
linear-programming method with input coming from other combinatorial methods can lead to improved result.








It is time to state the Frankl--R\"odl theorem.

\begin {theorem} [Frankl--R\"odl]

For every $\alpha,\beta,\gamma,\epsilon >0$, there is $\delta >0$ with the following property.
Let ${\cal U}_1$ be the family of $[\alpha n]$-subsets of $[n]$,
let ${\cal U}_2$ be the family of $[\beta n]$-subsets of $[n]$,
and let $X$ be the number of pairs of sets $A \in {\cal U}_1$, $B \in {\cal U}_2$
whose intersection is of size $[\gamma n]$.

Let $\cal F, \cal G$ be two subfamilies of ${\cal U}_1$ and ${\cal U}_2$, respectively,
with $|{\cal F}||{\cal G}| \ge (1-\delta)^n |{\cal U}_1|\cdot |{\cal U}_2|.$
Then the number of pairs $(A,B)$, $A \in {\cal F}$ and $B \in {\cal G}$ whose intersection has $[\gamma n]$ elements is at least
$(1-\epsilon)^nX$.

\end {theorem}

An important special case is where $\alpha=\beta=1/2$ and $\gamma=1/4$.
The Frankl--R\"odl paper contains generalizations in various directions. We could have assumed that,
e.g., $\cal F$ and $\cal G$ are  families of
partitions of $[n]$ into $r$ parts instead of families of sets. It also contains interesting
geometric applications. 
We will propose here two extensions of the Frankl--R\"odl theorem.

\subsection{Frankl--R\"odl/Frankl--Wilson with sum restrictions}

For $S \subset [n]$, we write $\|S\|=\sum\{s: s \in S\}$.
Let $\alpha_1,\alpha_2, \beta_1,\beta_2$ be reals such that
$0 < \alpha_1,\alpha_2 <1$, $0 < \beta_1,\beta_2 <1$. .
Consider the family $\cal G$ of subsets of $[n]$ such that
for every $S \in {\cal G}$ we have $|S|= [\alpha_1n]$, and $\|S\|=[\alpha_2 ({n \choose 2})]$.

Let $X$ be the number of pairs $A$ and $B$ in $\cal G$ with the properties:

(*) The intersection $C$ of $A$ and $B$ has precisely $[\beta_1 n]$ elements.

(**)  The sum of elements in $C$  is precisely $[\beta_2 ({n \choose 2})]$.

\begin {conj} [Frankl--R\"odl/Frankl--Wilson with sum restrictions]

For every $\epsilon >0$, there is $\delta >0$
such that if you have a subfamily
$F$ of $G$ of size $>(1-\delta)^n |G|$, then the number of pairs of sets in $F$
satisfying (*) and (**) is at least $(1-\epsilon)^n X$.

\end {conj}

\begin {rem}
(February 2015): Eoin Long has recently reduced many cases of this conjecture to the
original Frankl--R\"odl theorem.
\end {rem}

\subsection{Frankl--R\"odl/Frankl--Wilson for cocycles}

\begin {conj}[Frankl--R\"odl/Frankl--Wilson theorem for cocycles]

For every $\epsilon,\gamma >0$, there is $\delta >0$ with the following property.
Let $\cal F$ be the family of 3-cocycles. Let $X$ be
the number of pairs of elements in $\cal F$ whose symmetric difference has
precisely $m=[\gamma{{n} \choose {4}}]$
sets.
Then for every ${\cal G} \subset {\cal F}$ if $|{\cal G}| \ge (1-\delta)^{{n} \choose {4}} |{\cal F}|$, the
number of pairs of elements in $\cal G$ whose symmetric difference has
precisely $m$ sets is at least $(1-\epsilon)^{{n} \choose {4}} X$.

\end {conj}

The case of 1-cocycles is precisely the conclusion of Frankl--Wilson/Frankl--R\"odl
needed for Borsuk's conjecture, and a Frankl--R\"odl theorem for
4-cycles may also be a way to push up the asymptotic lower bounds for  Borsuk's
problem via Conjecture \ref {c:cycles}.



\section {Paul Erd\Horig{o}s' way with people and with mathematical problems }

There is a saying in the ancient Hebrew scriptures:

\begin {quotation}
Do not scorn any person and do not dismiss any thing, for there is no person
who has not his hour, and there is no thing that has not its place.
\end {quotation}

Paul Erd\Horig{o}s
had an amazing way of practicing this saying, when it came to people,
and likewise when it came to his beloved ``things,'' - mathematical problems.
And his way accounts for some of our finest hours.

\subsection *{Acknowledgment}
I am very thankful to
Evan DeCorte, Jeff Kahn, Andrey Kupavskii, Rom Pinchasi, Oleg Pikhurko, Andrei Raigorodskii, Micha Sharir, and Konrad Swanepoel for useful discussions and comments.

\myaddress

\newpage
\thispagestyle{empty}
\mbox{}
\newpage

\end{document}